# Classification and construction of picture fuzzy topological spaces

**Abdul Razaq[1,*], Harish Garg[2], Umer Shuaib[3]**

[1]Department of Mathematics, Division of Science and Technology, University of Education, Lahore 54770, Pakistan, abdul.razaq@ue.edu.pk

[2]School of Mathematics Thapar Institute of Engineering & Technology, Deemed University, Patiala 147004, Punjab, India; harish.garg@thapar.edu

[3]Department of Mathematics, Government College University, Faisalabad Pakistan; mumershuaib@gcuf.edu.pk

**∗ Correspondence:** abdul.razaq@ue.edu.pk

**Abstract.** This paper initiates the study of picture fuzzy topological spaces. In order to develop a mechanism to construct picture fuzzy topological spaces, we prove some basic results related to picture fuzzy sets together with the introduction of many new notions such as balanced and $\rho$-equivalent relations on picture fuzzy sets. We also define the rank, picture fuzzy base and picture fuzzy sub-base of picture fuzzy topological spaces. With the help of these notions, we present a method to design picture fuzzy topological spaces, which is the main outcome of this work. To the best of our knowledge, this is the first serious attempt to study picture fuzzy topological spaces. It is hoped that, this study will provide a decent platform for further development of picture fuzzy topological spaces and their applications in decision making theory.

## 1. Introduction

The development of crisp topology originated from classical analysis, and it has numerous applications in many research areas like data mining, machine learning, quantum gravity and data analysis [1–7]. The notion of topology discloses the relationship between spatial elements and features. It narrate some specific spatial mappings and designs datasets with superior quality control and better information integrity. In the recent years topology has not only firmly established itself as an important part of general mathematical structure but also of physicists mathematical arsenal. The notion the fuzzy topological space was introduced by Chang [8] in 1968. He generalized some fundamental concepts of topology such as open set, closed set, neighborhood and continuity. Based on this concept, many studies have been conducted in general theoretical areas and in different application sides. Lowen changed the fundamental characteristic of topology and presented one more definition of a fuzzy topological space [9,10]. Furthermore, in 1995, the concept of intuitionistic fuzzy topological space was introduced by Coker. He discussed counterparts of some general topological concepts like compactness and continuity [11]. For more details on intuitionistic fuzzy topological spaces, we suggest reading of [12–14]. Kramosil and Michálek presented a comprehensive study on fuzzy metrics and statistical metric spaces [15]. In [16,17], the of fuzzy soft topological space was introduced. The authors also developed a decision making method based on proposed fuzzy soft topological space.

This paper reveals the study of picture fuzzy topological spaces. The major contribution of this paper includes;

1. Some basic results regarding picture fuzzy inclusion, union and intersection have been discussed
2. Two novel relations on picture fuzzy subsets have been introduced. Based on these relations, a new notion, rank of picture fuzzy topological spaces, has been defined to classify picture fuzzy topological spaces.
3. The concepts of picture fuzzy base and picture fuzzy sub-base have been studied.

4. The final outcome of the (1), (2) and (3) is a mechanism to design various picture fuzzy topological spaces.

## 2. Fuzzy Picture Subsets

Zadeh [18] presented the concept of fuzzy set which has many applications on decision making, business, economy and data mining, etc. A fuzzy subset $\aleph$ on a universe $X$ is an object of the form $\{x, \mu_\aleph(x): x \in X\}$, where $\mu_\aleph: X \to [0,1]$ is called membership function of $\aleph$ and $\mu_\aleph(x)$ is known as degree of membership of $x$ in $\aleph$. Clearly, fuzzy set theory is the generalization of classical set theory. In classical sets, the corresponding membership function of a set $\aleph$ is the characteristic function of $\aleph$ that is equal to 1 $if$ $x \in \aleph$ and 0 if $x \notin \aleph$. After introduction of fuzzy sets, many new theories have been developed to deal with uncertainty and imprecision. Some of these theories are expansions of fuzzy set theory, while others endeavor to manage uncertainty and imprecision in some other suitable way. In 1983, Atanassov [19] presented the generalization fuzzy sets by defining intuitionistic fuzzy sets. An intuitionistic fuzzy subset $\aleph$ on a universe $X$ is an object of the form $\{x, \mu_\aleph(x), \rho_\aleph(x): x \in X\}$, where $\mu_\aleph: X \to [0,1]$ and $\rho_\aleph: X \to [0,1]$ are membership and non-membership functions respectively, such that $\mu_\aleph(x) + \rho_\aleph(x) \leq 1$ for all $x \in X$. In 2013, Cong [20] defined picture fuzzy sets which are one of the generalizations of intuitionistic fuzzy sets.

**Definition 1 [20]** Suppose $X$ is a non-empty crisp set and let us call $\mu_\aleph: X \to [0,1]$, $\rho_\aleph: X \to [0,1]$ and $\sigma_\aleph: X \to [0,1]$ positive, neutral and negative membership functions. A picture fuzzy subset $\aleph$ of $X$ is defined as $\aleph = \{x, \mu_\aleph(x), \rho_\aleph(x), \sigma_\aleph(x): x \in X\}$, where $\mu_\aleph(x) + \rho_\aleph(x) + \sigma_\aleph(x) \leq 1$. For all $x \in X$, $1 - (\mu_\aleph(x) + \rho_\aleph(x) + \sigma_\aleph(x))$ is termed as degree of refusal of $x$.

Picture fuzzy sets have variety of applications in decision-making environments. The development of picture fuzzy sets can be seen in [21]. In daily life, picture fuzzy set theory provides more than one choice for any decision. For instance,

**a.** Suppose a person is suffering from some disease. Then positive, neutral and negative membership functions can be associated with cureness, treatment and bitterness of disease respectively. Refusal can be related to the insufficient economic condition of the patient so that he cannot afford the expenses of hospital and refuses to be hospitalized.

**b.** Suppose a person has an allegation of a crime. Then positive, neutral and negative membership functions can be associated with maximum punishment, moderate punishment and release of the accused person respectively. Refusal can be related to the dismissal of the case due to reconciliation.

Let us rewrite the set operations over picture fuzzy subsets [20]:

**Definition 2** Let $\aleph_1$ and $\aleph_2$ be two picture fuzzy subsets of $X$, then

    i.    $\aleph_1 \subseteq \aleph_2$ if and only if $\mu_{\aleph_1}(x) \leq \mu_{\aleph_2}(x), \rho_{\aleph_1}(x) \leq \rho_{\aleph_2}(x)$ and $\sigma_{\aleph_1}(x) \geq \sigma_{\aleph_2}(x)$

ii. $\aleph_1 \uplus \aleph_2 = \{x, \max(\mu_{\aleph_1}(x), \mu_{\aleph_2}(x)), \min(\rho_{\aleph_1}(x), \rho_{\aleph_2}(x)), \min(\sigma_{\aleph_1}(x), \sigma_{\aleph_2}(x)) : x \in X\}$

iii. $\aleph_1 \cap \aleph_2 = \{x, \min(\mu_{\aleph_1}(x), \mu_{\aleph_2}(x)), \min(\rho_{\aleph_1}(x), \rho_{\aleph_2}(x)), \max(\sigma_{\aleph_1}(x), \sigma_{\aleph_2}(x)) : x \in X\}$

iv. $\aleph_1^c = \{x, \sigma_{\aleph_1}(x), \rho_{\aleph_1}(x), \mu_{\aleph_1}(x) : x \in X\}$

v. $\aleph_1 = \aleph_2$ if and only if $\aleph_1 \subseteq \aleph_2$ and $\aleph_2 \subseteq \aleph_1$

It should be noted that, in case of infinite union and intersection, maximum "max" and minimum "min" are replaced by supremum "sup" and infimum "inf" respectively.

Throughout this paper, we denote the family of all picture fuzzy subsets of $X$ by $\wp\mathcal{F}(X)$ and $I, O \in \wp\mathcal{F}(X)$ represent $\{x, 1, 0, 0 : x \in X\}$ and $\{x, 0, 0, 1 : x \in X\}$ respectively.

**Remark 1** Let $\aleph \in \wp\mathcal{F}(X)$. It is important to mention that $\aleph \subseteq I$ but $O \not\subseteq \aleph$ in general. Therefore, a question arises; for which family of picture fuzzy subsets $O \subseteq \aleph$ is true. We will give the answer later (see Theorem 11).

## 3. Picture fuzzy topology and some basic results for its construction

In this section, we define the notion of picture fuzzy topological space. Moreover, we present some basic results of picture fuzzy sets which are useful to design picture fuzzy topological spaces.

**Definition 3** Let $X$ denote a non-empty crisp set. A family $\mathfrak{I}$ of picture fuzzy subsets of $X$ is called a picture fuzzy topology on $X$ if

i. $I, O \in \mathfrak{I}$
ii. for any $\aleph_1, \aleph_2 \in \mathfrak{I}$, we have $\aleph_1 \cup \aleph_2 \in \mathfrak{I}$
iii. for any $\{\aleph_i\}_{i \in \Omega} \subseteq \mathfrak{I}$, we have $\bigcup_{i \in \Omega} \aleph_i \in \mathfrak{I}$, where $\Omega$ is an arbitrary index set.

If $\mathfrak{I}$ forms picture fuzzy topology on $X$, then the pair $(X, \mathfrak{I})$ denotes picture fuzzy topological space.

The following theorem is an immediate consequence of definition 2.

**Theorem 1.** For all picture fuzzy subsets $\aleph_1, \aleph_2$ and $\aleph_3$ of $X$

i. $\aleph_1 \uplus \aleph_2 = \aleph_2 \cap \aleph_1 = \aleph_1$ if $\aleph_1 = \aleph_2$
ii. $\aleph_1 \uplus \aleph_2 = \aleph_2 \uplus \aleph_1$ and $\aleph_1 \cap \aleph_2 = \aleph_2 \cap \aleph_1$
iii. $\aleph_1 \uplus (\aleph_2 \uplus \aleph_3) = (\aleph_1 \uplus \aleph_2) \uplus \aleph_3$ and $\aleph_1 \cap (\aleph_2 \cap \aleph_3) = (\aleph_1 \cap \aleph_2) \cap \aleph_3$

**Remark 2** The converse of Theorem 1 (i) is not true. Consider two picture fuzzy subsets $\aleph_1 = \{x, \mu_{\aleph_1}(x), \rho_{\aleph_1}(x), \sigma_{\aleph_1}(x) : x \in X\}$ and $\aleph_2 = \{x, \mu_{\aleph_2}(x), \rho_{\aleph_2}(x), \sigma_{\aleph_2}(x) : x \in X\}$ such that $\mu_{\aleph_1}(x) = \mu_{\aleph_2}(x), \rho_{\aleph_1}(x) \leq \rho_{\aleph_2}(x)$ and $\sigma_{\aleph_1}(x) = \sigma_{\aleph_2}(x)$. Then clearly $\aleph_1 \uplus \aleph_2 = \aleph_2 \cap \aleph_1 = \aleph_1$, but $\aleph_1 \neq \aleph_2$.

Remark 2 motivates us to define a new relation on $\wp\mathcal{F}(X)$.

**Definition 4** A picture fuzzy subset $\aleph_1$ is balanced with respect to picture fuzzy subset $\aleph_2$ if and only if $\mu_{\aleph_1}(x) = \mu_{\aleph_2}(x), \rho_{\aleph_1}(x) \leq \rho_{\aleph_2}(x)$ and $\sigma_{\aleph_1}(x) = \sigma_{\aleph_2}(x)$. If such is the case, then we write $\aleph_1 \equiv \aleph_2$.

The relation $\equiv$ is significant in designing various types of picture fuzzy topologies (see Theorem 13 and Example 8). Also, it is easy to verify that $\equiv$ is a partial order relation on $\wp\mathcal{F}(X)$.

The following result is the extension of Theorem 1 (i)

**Theorem 2** $\aleph_1 \uplus \aleph_2 = \aleph_2 \cap \aleph_1 = \aleph_1$ if and only if $\aleph_1 \equiv \aleph_2$

The proof is an immediate consequence of Definition 4.

The following Theorem shows that union and intersection are binary operations on $\wp F(X)$.

**Theorem 3** Let $\aleph_1 = \{x, \mu_{\aleph_1}(x), \rho_{\aleph_1}(x), \sigma_{\aleph_1}(x): x \in X\}$ and $\aleph_2 = \{x, \mu_{\aleph_2}(x), \rho_{\aleph_2}(x), \sigma_{\aleph_2}(x): x \in X\}$ be two fuzzy subsets of $X$. Then $\aleph_1 \cup \aleph_2$ and $\aleph_1 \cap \aleph_2$ are picture fuzzy subsets of $X$.

**Proof.** We know

$$\aleph_1 \cup \aleph_2 = \{x, \mu_{\aleph_1 \cup \aleph_2}(x), \rho_{\aleph_1 \cup \aleph_2}(x), \sigma_{\aleph_1 \cup \aleph_2}(x): x \in X\}$$

$$= \{x, \max(\mu_{\aleph_1}(x), \mu_{\aleph_2}(x)), \min(\rho_{\aleph_1}(x), \rho_{\aleph_2}(x)), \min(\sigma_{\aleph_1}(x), \sigma_{\aleph_2}(x)): x \in X\} \text{ and}$$

$$\aleph_1 \cap \aleph_2 = \{x, \mu_{\aleph_1 \cap \aleph_2}(x), \rho_{\aleph_1 \cap \aleph_2}(x), \sigma_{\aleph_1 \cap \aleph_2}(x): x \in X\}$$

$$= \{x, \min(\mu_{\aleph_1}(x), \mu_{\aleph_2}(x)), \min(\rho_{\aleph_1}(x), \rho_{\aleph_2}(x)), \max(\sigma_{\aleph_1}(x), \sigma_{\aleph_2}(x)): x \in X\}$$

For all pairs $(\mu_{\aleph_1}(x), \mu_{\aleph_2}(x))$, $(\rho_{\aleph_1}(x), \rho_{\aleph_2}(x))$ and $(\sigma_{\aleph_1}(x), \sigma_{\aleph_2}(x))$, there are following eight possibilities;

i. $\mu_{\aleph_1}(x) \geq \mu_{\aleph_2}(x), \rho_{\aleph_1}(x) \geq \rho_{\aleph_2}(x)$ and $\sigma_{\aleph_1}(x) \geq \sigma_{\aleph_2}(x)$
ii. $\mu_{\aleph_1}(x) \geq \mu_{\aleph_2}(x), \rho_{\aleph_1}(x) \leq \rho_{\aleph_2}(x)$ and $\sigma_{\aleph_1}(x) \geq \sigma_{\aleph_2}(x)$
iii. $\mu_{\aleph_1}(x) \geq \mu_{\aleph_2}(x), \rho_{\aleph_1}(x) \geq \rho_{\aleph_2}(x)$ and $\sigma_{\aleph_1}(x) \leq \sigma_{\aleph_2}(x)$
iv. $\mu_{\aleph_1}(x) \geq \mu_{\aleph_2}(x), \rho_{\aleph_1}(x) \leq \rho_{\aleph_2}(x)$ and $\sigma_{\aleph_1}(x) \leq \sigma_{\aleph_2}(x)$
v. $\mu_{\aleph_1}(x) \leq \mu_{\aleph_2}(x), \rho_{\aleph_1}(x) \leq \rho_{\aleph_2}(x)$ and $\sigma_{\aleph_1}(x) \leq \sigma_{\aleph_2}(x)$
vi. $\mu_{\aleph_1}(x) \leq \mu_{\aleph_2}(x), \rho_{\aleph_1}(x) \geq \rho_{\aleph_2}(x)$ and $\sigma_{\aleph_1}(x) \leq \sigma_{\aleph_2}(x)$
vii. $\mu_{\aleph_1}(x) \leq \mu_{\aleph_2}(x), \rho_{\aleph_1}(x) \leq \rho_{\aleph_2}(x)$ and $\sigma_{\aleph_1}(x) \geq \sigma_{\aleph_2}(x)$
viii. $\mu_{\aleph_1}(x) \leq \mu_{\aleph_2}(x), \rho_{\aleph_1}(x) \geq \rho_{\aleph_2}(x)$ and $\sigma_{\aleph_1}(x) \geq \sigma_{\aleph_2}(x)$

It can be easily verified that $\mu_{\aleph_1 \cup \aleph_2} + \rho_{\aleph_1 \cup \aleph_2} + \sigma_{A \cup B} \leq 1$ and $\mu_{\aleph_1 \cap \aleph_2} + \rho_{\aleph_1 \cap \aleph_2} + \sigma_{\aleph_1 \cap \aleph_2} \leq 1$ for all above possibilities.

The following theorem shows that distributive laws for union and intersection holds in $\wp F(X)$.

**Theorem 4** For all picture fuzzy subsets $\aleph_1, \aleph_2$ and $\aleph_3$ of $X$

i. $\aleph_1 \cup (\aleph_2 \cap \aleph_3) = (\aleph_1 \cup \aleph_2) \cap (\aleph_1 \cup \aleph_3)$
ii. $\aleph_1 \cap (\aleph_2 \cup \aleph_3) = (\aleph_1 \cap \aleph_2) \cup (\aleph_1 \cap \aleph_3)$

**Prrof.** i. Let $x \in X$, then

$$\mu_{\aleph_1 \cup (\aleph_2 \cap \aleph_3)}(x) = \max(\mu_{\aleph_1}(x), \mu_{\aleph_2 \cap \aleph_3}(x)) = \max(\mu_{\aleph_1}(x), \min(\mu_{\aleph_2}(x), \mu_{\aleph_3}(x))) \quad (3.1)$$

and

$$\mu_{(\aleph_1 \cup \aleph_2) \cap (\aleph_1 \cup \aleph_3)}(x) = \min(\mu_{\aleph_1 \cup \aleph_2}(x), \mu_{\aleph_1 \cup \aleph_3}(x)) =$$
$$\min(\max(\mu_{\aleph_1}(x), \mu_{\aleph_2}(x)), \max(\mu_{\aleph_1}(x), \mu_{\aleph_3}(x))) \quad (3.2)$$

Now there are three possibilities

a) If $max\left(\mu_{\aleph_1}(x), \mu_{\aleph_2}(x), \mu_{\aleph_3}(x)\right) = \mu_{\aleph_1}(x)$, then (3.1) and (3.2) yield

$$\mu_{\aleph_1 \cup (\aleph_2 \cap \aleph_3)}(x) = \mu_{\aleph_1}(x) = \mu_{(\aleph_1 \cup \aleph_2) \cap (\aleph_1 \cup \aleph_3)}(x)$$

b) If $max\left(\mu_{\aleph_1}(x), \mu_{\aleph_2}(x), \mu_{\aleph_3}(x)\right) = \mu_{\aleph_2}(x)$, then from (1) and (2), we obtain

$$\mu_{\aleph_1 \cup (\aleph_2 \cap \aleph_3)}(x) = max\left(\mu_{\aleph_1}(x), \mu_{\aleph_3}(x)\right) = \mu_{(\aleph_1 \cup \aleph_2) \cap (\aleph_1 \cup \aleph_3)}(x)$$

c) If $max\left(\mu_{\aleph_1}(x), \mu_{\aleph_2}(x), \mu_{\aleph_3}(x)\right) = \mu_{\aleph_3}(x)$, then (1) and (2) imply that

$$\mu_{\aleph_1 \cup (\aleph_2 \cap \aleph_3)}(x) = max\left(\mu_{\aleph_1}(x), \mu_{\aleph_2}(x)\right) = \mu_{(\aleph_1 \cup \aleph_2) \cap (\aleph_1 \cup \aleph_3)}(x)$$

Thus, in all three cases we have $\mu_{\aleph_1 \cup (\aleph_2 \cap \aleph_3)}(x) = \mu_{(\aleph_1 \cup \aleph_2) \cap (\aleph_1 \cup \aleph_3)}(x)$. Similarly, we can show $\rho_{\aleph_1 \cup (\aleph_2 \cap \aleph_3)}(x) = \rho_{(\aleph_1 \cup \aleph_2) \cap (\aleph_1 \cup \aleph_3)}(x)$ and $\sigma_{\aleph_1 \cup (\aleph_2 \cap \aleph_3)}(x) = \sigma_{(\aleph_1 \cup \aleph_2) \cap (\aleph_1 \cup \aleph_3)}(x)$. Hence $\aleph_1 \cup (\aleph_2 \cap \aleph_3) = (\aleph_1 \cup \aleph_2) \cap (\aleph_1 \cup \aleph_3)$.

ii. The proof is similar to that of (i)"

## 4. Equivalence classes in $\wp\mathcal{F}(X)$ and rank of picture fuzzy topology

In this section, we define an important equivalence relation $\|$ on $\wp\mathcal{F}(X)$. We show that all picture fuzzy subsets belong to the same equivalence class under $\|$ satisfy some basic properties of inclusion, intersection and union (See Theorem 11). We also classify picture fuzzy topological spaces into different ranks on the basis of newly defined equivalence relation $\|$.

We begin this section with the following theorem.

**Theorem 5** Let $\aleph_1 = \{x, \mu_{\aleph_1}(x), \rho_{\aleph_1}(x), \sigma_{\aleph_1}(x) : x \in X\}$ and $\aleph_2 = \{x, \mu_{\aleph_2}(x), \rho_{\aleph_2}(x), \sigma_{\aleph_2}(x) : x \in X\}$ be two fuzzy subsets of $X$. Then

i. $\aleph_1 \cap \aleph_2 \subseteq \aleph_1$ and $\aleph_1 \cap \aleph_2 \subseteq \aleph_2$      (4.1)
ii. $\aleph_1 \subseteq \aleph_2$ if and only if $\aleph_1 \cap \aleph_2 = \aleph_1$      (4.2)

**Proof.** We know

$$\aleph_1 \uplus \aleph_2 = \{x, \mu_{\aleph_1 \cup \aleph_2}(x), \rho_{\aleph_1 \cup \aleph_2}(x), \sigma_{\aleph_1 \cup \aleph_2}(x) : x \in X\}$$

$$= \{x, max\left(\mu_{\aleph_1}(x), \mu_{\aleph_2}(x)\right), min\left(\rho_{\aleph_1}(x), \rho_{\aleph_2}(x)\right), min\left(\sigma_{\aleph_1}(x), \sigma_{\aleph_2}(x)\right) : x \in X\}$$ and

$$\aleph_1 \cap \aleph_2 = \{x, \mu_{\aleph_1 \cap \aleph_2}(x), \rho_{\aleph_1 \cap \aleph_2}(x), \sigma_{\aleph_1 \cap \aleph_2}(x) : x \in X\}$$

$$= \{x, min\left(\mu_{\aleph_1}(x), \mu_{\aleph_2}(x)\right), min\left(\rho_{\aleph_1}(x), \rho_{\aleph_2}(x)\right), max\left(\sigma_{\aleph_1}(x), \sigma_{\aleph_2}(x)\right) : x \in X\}$$

i) Since $\mu_{\aleph_1 \cap \aleph_2}(x) = min\left(\mu_{\aleph_1}(x), \mu_{\aleph_2}(x)\right)$, $\rho_{\aleph_1 \cap \aleph_2}(x) = min\left(\rho_{\aleph_1}(x), \rho_{\aleph_2}(x)\right)$ and $\sigma_{\aleph_1 \cap \aleph_2}(x) = max\left(\sigma_{\aleph_1}(x), \sigma_{\aleph_2}(x)\right)$. This means that $\mu_{\aleph_1 \cap \aleph_2}(x) \leq \mu_{\aleph_1}(x), \mu_{\aleph_2}(x)$, $\rho_{\aleph_1 \cap \aleph_2}(x) \leq \rho_{\aleph_1}(x), \rho_{\aleph_2}(x)$ and $\sigma_{\aleph_1 \cap \aleph_2}(x) \geq \sigma_{\aleph_1}(x), \sigma_{\aleph_2}(x)$. Thus, $\aleph_1 \cap \aleph_2 \subseteq \aleph_1$ and $\aleph_1 \cap \aleph_2 \subseteq \aleph_2$.

ii) Let $\aleph_1 \subseteq \aleph_2$, then $\mu_{\aleph_1}(x) \leq \mu_{\aleph_2}(x)$, $\rho_{\aleph_1}(x) \leq \rho_{\aleph_2}(x)$ and $\sigma_{\aleph_1}(x) \geq \sigma_{\aleph_2}(x)$. It reveals that

$$\mu_{\aleph_1 \cap \aleph_2}(x) = \mu_{\aleph_1}(x), \rho_{\aleph_1 \cap \aleph_2}(x) = \rho_{\aleph_1}(x) \text{ and } \sigma_{\aleph_1 \cap \aleph_2}(x) = \sigma_{\aleph_1}(x)$$

implying that $\aleph_1 \cap \aleph_2 = \aleph_1$.

Conversely, suppose that $\aleph_1 \cap \aleph_2 = \aleph_1$. Since $\aleph_1 \cap \aleph_2 \subseteq \aleph_2$, therefore $\aleph_1 \subseteq \aleph_2$.

The following examples show that it is not necessary for picture fuzzy subsets $\aleph_1$ and $\aleph_2$ such that

i. $\aleph_1 \subseteq \aleph_1 \uplus \aleph_2$ and $\aleph_2 \subseteq \aleph_1 \uplus \aleph_2$ (4.3)

ii. $\aleph_1 \subseteq \aleph_2$ if and only if $\aleph_1 \uplus \aleph_2 = \aleph_2$ (4.4)

**Example 1** Let $X = \{a, b, c\}$ and $\aleph_1 = \{(a, 0.50, 0.20, 0.25), (b, 0.40, 0.10, 0.50), (c, 0.20, 0.30, 0.45)\}$, $\aleph_2 = \{(a, 0.40, 0.30, 0.10), (b, 0.20, 0.60, 0.10), (c, 0.30, 0.20, 0.15)\}$ be two picture fuzzy subsets of $X$. Then $\aleph_1 \uplus \aleph_2 = \{(a, 0.50, 0.20, 0.10), (b, 0.40, 0.10, 0.10), (c, 0.30, 0.20, 0.15)\}$. Clearly, neither $\aleph_1 \subseteq \aleph_1 \uplus \aleph_2$ nor $\aleph_2 \subseteq \aleph_1 \uplus \aleph_2$.

**Example 2** Consider two picture fuzzy subsets $\aleph_1 = \{(a, 0.30, 0.20, 0.25), (b, 0.10, 0.30, 0.50), (c, 0.20, 0.20, 0.45)\}$ and $\aleph_2 = \{(a, 0.40, 0.30, 0.10), (b, 0.20, 0.60, 0.10), (c, 0.30, 0.20, 0.15)\}$ of $X = \{a, b, c\}$. Then $\aleph_1 \uplus \aleph_2 = \{(a, 0.40, 0.20, 0.10), (b, 0.20, 0.30, 0.10), (c, 0.30, 0.20, 0.15)\}$. Clearly, $\aleph_1 \uplus \aleph_2 \neq \aleph_2$.

**Example 3** Consider $\aleph_1 = \{(a, 0.30, 0.20, 0.25), (b, 0.10, 0.30, 0.50), (c, 0.20, 0.20, 0.45)\}$ and $\aleph_2 = \{(a, 0.40, 0.15, 0.10), (b, 0.20, 0.25, 0.10), (c, 0.30, 0.20, 0.15)\}$. Then $\aleph_1 \uplus \aleph_2 = \aleph_2$ but $\aleph_1 \nsubseteq \aleph_2$.

Thus, a curiosity arises to explore the condition for the validation of the statement (4.3) and (4.4). The upcoming three theorems resolve this curiosity.

**Theorem 7** $\aleph_1 \subseteq \aleph_1 \uplus \aleph_2$ if and only if $\rho_{\aleph_1}(x) \leq \rho_{\aleph_2}(x)$

**Proof** Let $\aleph_1 \subseteq \aleph_1 \uplus \aleph_2$, then $\rho_{\aleph_1}(x) \leq \rho_{\aleph_1 \uplus \aleph_2}(x) = min\left(\rho_{\aleph_1}(x), \rho_{\aleph_2}(x)\right)$. Thus, $\rho_{\aleph_1}(x) \leq \rho_{\aleph_2}(x)$.

Conversely, suppose that $\rho_{\aleph_1}(x) \leq \rho_{\aleph_2}(x)$. Then $\mu_{\aleph_1 \uplus \aleph_2}(x) = max\left(\mu_{\aleph_1}(x), \mu_{\aleph_2}(x)\right) \geq \mu_{\aleph_1}(x)$, $\rho_{\aleph_1 \uplus \aleph_2}(x) = min\left(\rho_{\aleph_1}(x), \rho_{\aleph_2}(x)\right) = \rho_{\aleph_1}(x)$ and $\sigma_{\aleph_1 \uplus \aleph_2}(x) = min\left(\sigma_{\aleph_1}(x), \sigma_{\aleph_2}(x)\right) \leq \sigma_{\aleph_1}(x)$, Thus, $\aleph_1 \subseteq \aleph_1 \uplus \aleph_2$.

**Theorem 8** $\aleph_1 \subseteq \aleph_1 \uplus \aleph_2$ and $\aleph_2 \subseteq \aleph_1 \uplus \aleph_2$ and if and only if $\rho_{\aleph_1}(x) = \rho_{\aleph_2}(x)$

**Proof.** The proof is an immediate consequence of Theorem 7.

**Theorem 9** Let $\rho_{\aleph_1}(x) = \rho_{\aleph_2}(x)$ for all $x \in X$. Then $\aleph_1 \subseteq \aleph_2$ if and only if $\aleph_1 \uplus \aleph_2 = \aleph_2$.

**Proof.** Let $\rho_{\aleph_1}(x) = \rho_{\aleph_2}(x)$ for all $x \in X$. Then $\aleph_1 \subseteq \aleph_2 \Leftrightarrow max\left(\mu_{\aleph_1}(x), \mu_{\aleph_2}(x)\right) = \mu_{\aleph_2}(x)$, $min\left(\rho_{\aleph_1}(x), \rho_{\aleph_2}(x)\right) = \rho_{\aleph_1}(x) = \rho_{\aleph_2}(x)$ and $min\left(\sigma_{\aleph_1}(x), \sigma_{\aleph_2}(x)\right) = \sigma_{\aleph_2}(x) \Leftrightarrow \aleph_1 \uplus \aleph_2 = \aleph_2$

The theorems 8 and 9, reveal that the condition for picture fuzzy subsets $\aleph_1$ and $\aleph_2$ to satisfy conditions (4.3) and (4.4) is $\rho_{\aleph_1}(x) = \rho_{\aleph_2}(x)$. This leads us to define a new relation on $\wp\mathcal{F}(X)$.

**Definition 5** A picture fuzzy subset $\aleph_1 = \{x, \mu_{\aleph_1}(x), \rho_{\aleph_1}(x), \sigma_{\aleph_1}(x) : x \in X\}$ is $\rho$-equivalent to a picture fuzzy subset $\aleph_2 = \{x, \mu_{\aleph_2}(x), \rho_{\aleph_2}(x), \sigma_{\aleph_2}(x) : x \in X\}$ if $\mu_{\aleph_1}(x) = \rho_{\aleph_2}(x)$ for all $x \in X$. If $\aleph_1$ is $\rho$-equivalent to $\aleph_2$, then we write $\aleph_1 \parallel \aleph_2$.

It is easy to verify that ∥ is an equivalence relation and hence all equivalences classes under ∥ form a partition of the set $\wp\mathcal{F}(X)$ all picture fuzzy subsets of $X$.

**Theorem 10** $\aleph_1 = \aleph_2$ if and only if $\aleph_1 \parallel \aleph_2$ and $\aleph_1 \equiv \aleph_2$

The proof is straightforward

The final conclusion of this section can be summarized in the following theorem

**Theorem 11** $\aleph_1$ and $\aleph_2$ belong to the same equivalence class under ∥, if and only if the followings are true.

i. $\aleph_1 \cap \aleph_2 \subseteq \aleph_1$ and $\aleph_1 \cap \aleph_2 \subseteq \aleph_2$
ii. $\aleph_1 \subseteq \aleph_2$ if and only if $\aleph_1 \cap \aleph_2 = \aleph_1$.
iii. $\aleph_1 \subseteq \aleph_1 \cup \aleph_2$ and $\aleph_2 \subseteq \aleph_1 \cup \aleph_2$
iv. $\aleph_1 \subseteq \aleph_2$ if and only if $\aleph_1 \cup \aleph_2 = \aleph_2$
v. $\aleph_1 \subseteq I$ and $O \subseteq \aleph_1$

The Theorem **11** plays a vital role in the construction of picture topological spaces. Since it holds for the picture fuzzy subsets of the same equivalence class, therefore, we classify picture fuzzy topologies on the basis of the number equivalence classes. This leads us to define an important notion called rank of picture fuzzy topology.

**Definition 6** A natural number $n$ is called rank of picture fuzzy topology $\mathfrak{I}$ if it has $n$ distinct equivalence classes under ∥.

## 5. Picture fuzzy bases and sub-bases and construction of picture fuzzy topological spaces

In this section, firstly we define the notions of picture fuzzy bases and sub-bases of picture fuzzy topology. Next, we construct some picture fuzzy topologies of different picture fuzzy sub-bases and ranks.

**Definition 7** Let $\mathcal{B} = \{\beta_i \in \wp\mathcal{F}(X)\}$ such that $I, O \notin \mathcal{B}$. Consider another sub-collection $\mathfrak{I}$ of $\wp\mathcal{F}(X)$ with the following properties;

i. $I, O \in \mathfrak{I}$

ii. $\mathcal{V}_i \in \mathfrak{I} \setminus \{O, I\}$ and $\mathcal{V}_i \nparallel O$ if and only if there exist $\beta_i \in \mathcal{B}$ such that if $\mathcal{V}_i = \cup_{i \in \Omega} \beta_i$.

iii. $\mathcal{V}_i \in \mathfrak{I} \setminus \{O, I\}$ and $\mathcal{V}_i \parallel O$ if and only there exist $\beta_i \in \mathcal{B}$ such that $\mathcal{V}_i = O \cup \{\cup_{i \in \Omega} \beta_i\}$.

If $\mathfrak{I}$ forms picture fuzzy topology on $X$, then $\mathcal{B}$ is called a picture fuzzy base for $\mathfrak{I}$.

**Theorem 12** $\mathcal{B}$ is a picture fuzzy base for $\mathfrak{I}$ if it is closed under finite intersection.

**Proof** Suppose that $\mathcal{B}$ is closed under finite intersection. We want to show that $\mathfrak{I}$ is a picture fuzzy topology on $X$.

(1) Definition 7 shows that $I, O \in \mathfrak{I}$.

(2) $\mathcal{V}_1, \mathcal{V}_2 \in \mathfrak{I} \setminus \{O, I\}$, then there are three cases (i) $\mathcal{V}_1, \mathcal{V}_2 \nparallel$ (ii) $\mathcal{V}_1 \parallel O$ and $\mathcal{V}_2 \nparallel O$ (iii) $\mathcal{V}_1, \mathcal{V}_2 \parallel O$

Let (i) be true, then $\mathcal{V}_1 \cap \mathcal{V}_2 = (\beta_1 \cup \beta_2 \cup \ldots \cup \beta_r \ldots) \cap (\beta'_1 \cup \beta'_2 \cup \ldots \cup \beta'_s \ldots)$

$= (\beta_1 \cap \beta_1') \cup (\beta_1 \cap \beta_2') \cup \ldots \cup (\beta_1 \cap \beta_s') \cup \ldots \cup (\beta_2 \cap \beta_1') \cup (\beta_2 \cap \beta_2') \cup \ldots \cup (\beta_2 \cap \beta_s') \cup \ldots \cup (\beta_r \cap \beta_1') \cup (\beta_r \cap \beta_2') \cup \ldots \cup (\beta_r \cap \beta_s') \cup \ldots = \cup_{i,j \in \Omega}(\beta_i \cap \beta_j')$, where $\beta_i, \beta_j' \in \mathcal{B}$.

Since $\mathcal{B}$ is closed under finite intersection, therefore each $\beta_j \cap \beta_j' \in \mathcal{B}$. Thus, $\mathcal{V}_1 \cap \mathcal{V}_2 = \cup_{i,j \in \Omega}(\beta_i \cap \beta_j') \in \mathfrak{J}$.

Similarly, $\mathcal{V}_1 \cap \mathcal{V}_2 \in \mathfrak{J}$ for the other (ii) and (iii) as well.

(3) Suppose $\{\mathcal{V}_i\}_{i \in \Omega} \subseteq \mathfrak{J} \setminus \{O, I\}$ and if
(i) each $\mathcal{V}_i \not\equiv O$, then $\cup_{i \in \Omega} \mathcal{V}_i = \cup_{j \in \Omega'} \beta_j \in \mathfrak{J}$

(ii) some of the $\mathcal{V}_k \in \{\mathcal{V}_i\}_{i \in \Omega}$ are $\rho$-equivalent to $O$, then $\cup_{i \in \Omega} \mathcal{V}_i = O \cup \left(\cup_{j \in \Omega'} \beta_j\right) \in \mathfrak{J}$

Next, $O \cap \mathcal{V}_i = I \cup \mathcal{V}_i$ and $O \cup \mathcal{V}_i = \mathcal{V}_i = O \cup \mathcal{V}_i$ together with (2) and (3) imply;

for all $\mathcal{V}_1, \mathcal{V}_2 \in \mathfrak{J}$, we have $\mathcal{V}_1 \cap \mathcal{V}_2 \in \mathfrak{J}$ and for all $\{\mathcal{V}_i\}_{i \in \Omega} \subseteq \mathfrak{J}$, we have $\cup_{i \in \Omega} \mathcal{V}_i \in \mathfrak{J}$.

**Definition 8** A collection $\mathfrak{S}$ of picture fuzzy subset of $X$ is called picture fuzzy sub-base of picture fuzzy topology $\mathfrak{J}$ if

i. all finite intersections of the elements of $\mathfrak{S}$ evolves a base for some picture fuzzy topology $\mathfrak{J}$.

ii. $\aleph_1 \cap \aleph_2 \in \mathfrak{S}$ implies $\aleph_1 \subseteq \aleph_2$ for all $\aleph_1, \aleph_2 \in \mathfrak{S}$

It should be that, we want to see minimum elements in sub-base that is why we impose condition (ii) in its definition.

**Remark 3** Every collection $\mathfrak{S} = \{\aleph_i : \aleph_i \in \wp\mathcal{F}(X)\}$ of picture fuzzy subset of $X$, satisfying condition (ii) of Definition 8, yields some picture fuzzy topology $\mathfrak{J}$. Because the collection $\mathcal{B}$, evolved by all finite intersections of the elements of $\mathfrak{S}$, is closed under finite intersection and hence Theorem 12 guarantees its role as a base for picture fuzzy topology $\mathfrak{J}$.

**Remark 4** The arbitrary union and finite intersection of picture fuzzy subsets from one equivalence class belong to the same equivalence class. This together with definitions of picture fuzzy base and picture fuzzy sub-base implies that "the sub-base of picture fuzzy topology $\mathfrak{J}$ of rank $n$ contains at least $n - 1$ elements.

### 5.1 Construction of picture fuzzy topological spaces

All the elements in picture fuzzy topologies of rank 1 belong to the same equivalence class. Therefore, Theorem 12 is valid for these elements. For this reason, we found that it is relatively easy to construct the picture fuzzy topologies of rank 1 as compared to those of higher ranks.

The following facts regarding in picture fuzzy topology $\mathfrak{J}$ of rank 1 are easy to prove;

i. $\{I = \{x, 1,0,0\}, O = \{x, 0,0,1\}\}$ is the smallest picture fuzzy topology of type 1. It is a trivial picture fuzzy topology.

ii. The picture fuzzy topology of rank 1 with sub-base $\mathfrak{S} = \{\aleph_1\}$ is $\{I, O, \aleph_1\}$.

iii. The picture fuzzy topology of rank 1 with $\mathfrak{S} = \{\aleph_1, \aleph_2, \ldots, \aleph_n\}$ such that $\aleph_1 \subset \aleph_2 \subset \cdots \subset \aleph_n$ is $\{I, O, \aleph_1, \aleph_2, \ldots, \aleph_n\}$.

The upcoming example shows the structure of picture fuzzy topology of rank 2 on $X = \{a, b, c\}$ with sub-base $\mathfrak{S} = \{\aleph_1, \aleph_2\}$ such that neither $\aleph_1 \subset \aleph_2$ nor $\aleph_2 \subset \aleph_1$.

**Example 4** Consider $\mathfrak{S} = \begin{cases} \aleph_1 = \{(a, 0.25, 0.20, 0.30), (b, 0.35, 0.10, 0.45), (c, 0.30, 0.35, 0.10)\}, \\ \aleph_2 = \{(a, 0.45, 0.20, 0.35), (b, 0.25, 0.10, 0.40), (c, 0.50, 0.35, 0.05)\} \end{cases}$. Then all finite intersections of the elements of $\mathfrak{S}$ give $\mathcal{B} =$
$\begin{cases} \aleph_1 = \{(a, 0.25, 0.20, 0.30), (b, 0.35, 0.10, 0.45), (c, 0.30, 0.35, 0.10)\}, \\ \aleph_2 = \{(a, 0.45, 0.20, 0.35), (b, 0.25, 0.10, 0.40), (c, 0.50, 0.35, 0.05)\}, \\ \aleph_1 \cap \aleph_2 = \\ \aleph_3 = \{(a, 0.25, 0.20, 0.35), (b, 0.25, 0.10, 0.45), (c, 0.30, 0.35, 0.10)\} \end{cases}$. The picture fuzzy topology obtained from $\mathcal{B}$ is

$$\mathfrak{I} = \begin{cases} I = \{(a, 1.00, 0.00, 0.00), (b, 1.00, 0.00, 0.00), (c, 1.00, 0.00, 0.00)\}, \\ O = \{(a, 0.00, 0.00, 1.00), (b, 0.00, 0.00, 1.00), (c, 0.00, 0.00, 1.00)\} \\ \aleph_1 = \{(a, 0.25, 0.20, 0.30), (b, 0.35, 0.10, 0.45), (c, 0.30, 0.35, 0.10)\}, \\ \aleph_2 = \{(a, 0.45, 0.20, 0.35), (b, 0.25, 0.10, 0.40), (c, 0.50, 0.35, 0.05)\}, \\ \aleph_1 \cap \aleph_2 = \\ \aleph_3 = \{(a, 0.25, 0.20, 0.35), (b, 0.25, 0.10, 0.45), (c, 0.30, 0.35, 0.10)\}, \\ \aleph_1 \cup \aleph_2 = \\ \aleph_4 = \{(a, 0.45, 0.20, 0.30), (b, 0.35, 0.10, 0.40), (c, 0.50, 0.35, 0.05)\}, \\ O \cup \aleph_1 = \\ \aleph_5 = \{(a, 0.25, 0.00, 0.30), (b, 0.35, 0.00, 0.45), (c, 0.30, 0.00, 0.10)\}, \\ O \cup \aleph_2 = \\ \aleph_6 = \{(a, 0.45, 0.00, 0.35), (b, 0.25, 0.00, 0.40), (c, 0.50, 0.00, 0.05)\}, \\ O \cup \aleph_3 = \\ \aleph_7 = \{(a, 0.25, 0.00, 0.35), (b, 0.25, 0.00, 0.45), (c, 0.30, 0.00, 0.10)\}, \\ O \cup \aleph_4 = \\ \aleph_8 = \{(a, 0.45, 0.00, 0.30), (b, 0.35, 0.00, 0.40), (c, 0.50, 0.00, 0.05)\} \end{cases}.$$

In the following example, we construct the picture fuzzy topology of rank 3 on $X = \{a, b, c\}$ with sub-base $\mathfrak{S} = \{\aleph_1, \aleph_2\}$ such that $\aleph_1 \subset \aleph_2$ but $\aleph_1 \not\equiv \aleph_2$.

**Example 5** Consider $\mathfrak{S} = \begin{cases} \aleph_1 = \{(a, 0.30, 0.20, 0.45), (b, 0.20, 0.25, 0.40), (c, 0.30, 0.35, 0.10)\}, \\ \aleph_2 = \{(a, 0.45, 0.30, 0.35), (b, 0.25, 0.30, 0.30), (c, 0.50, 0.40, 0.05)\} \end{cases}$. The base $\mathcal{B}$ evolved from $\mathfrak{S}$ is $\mathcal{B}$ itself. Thus, picture fuzzy topology constructed from $\mathcal{B}$ is

$$\mathfrak{I} = \begin{cases} I = \{(a, 1.00, 0.00, 0.00), (b, 1.00, 0.00, 0.00), (c, 1.00, 0.00, 0.00)\}, \\ O = \{(a, 0.00, 0.00, 1.00), (b, 0.00, 0.00, 1.00), (c, 0.00, 0.00, 1.00)\} \\ \aleph_1 = \{(a, 0.30, 0.20, 0.45), (b, 0.20, 0.25, 0.40), (c, 0.30, 0.35, 0.10)\}, \\ \aleph_2 = \{(a, 0.45, 0.30, 0.35), (b, 0.25, 0.30, 0.30), (c, 0.50, 0.40, 0.05)\}, \\ \aleph_1 \cup \aleph_2 = \\ \aleph_3 = \{(a, 0.45, 0.20, 0.35), (b, 0.25, 0.25, 0.30), (c, 0.50, 0.35, 0.05)\}, \\ O \cup \aleph_1 = \\ \aleph_4 = \{(a, 0.30, 0.00, 0.45), (b, 0.20, 0.00, 0.40), (c, 0.30, 0.00, 0.10)\}, \\ O \cup \aleph_2 = \\ \aleph_5 = \{(a, 0.45, 0.00, 0.35), (b, 0.25, 0.00, 0.30), (c, 0.50, 0.00, 0.05)\} \end{cases}.$$

The succeeding example describes the design of picture fuzzy topology of rank 3 on $X = \{a, b, c\}$ with sub-base $\mathfrak{S} = \{\aleph_1, \aleph_2\}$ such that $\aleph_1 \subset \aleph_2$ and $\aleph_1 \equiv \aleph_2$.

**Example 6** Consider $= \begin{cases} \aleph_1 = \{(a, 0.35, 0.20, 0.25), (b, 0.20, 0.15, 0.30), (c, 0.20, 0.35, 0.15)\}, \\ \aleph_2 = \{(a, 0.35, 0.30, 0.25), (b, 0.20, 0.25, 0.30), (c, 0.20, 0.40, 0.15)\} \end{cases}$. The base $\mathcal{B}$ evolved from $\mathfrak{S}$ is $\mathcal{B}$ itself. We obtain the following picture fuzzy topology from $\mathcal{B}$;

$$\mathfrak{I} = \begin{cases} I = \{(a, 1.00, 0.00, 0.00), (b, 1.00, 0.00, 0.00), (c, 1.00, 0.00, 0.00)\}, \\ O = \{(a, 0.00, 0.00, 1.00), (b, 0.00, 0.00, 1.00), (c, 0.00, 0.00, 1.00)\} \\ \aleph_1 = \{(a, 0.35, 0.20, 0.25), (b, 0.20, 0.15, 0.30), (c, 0.20, 0.35, 0.15)\}, \\ \aleph_2 = \{(a, 0.35, 0.30, 0.25), (b, 0.20, 0.25, 0.30), (c, 0.20, 0.40, 0.15)\}, \\ O \cup \aleph_1 = \\ \aleph_3 = \{(a, 0.35, 0.00, 0.25), (b, 0.20, 0.00, 0.30), (c, 0.20, 0.00, 0.15)\} \end{cases}.$$

Next, we form a picture fuzzy topology of rank 3 on $X = \{a, b, c\}$ with sub-base $\mathfrak{S} = \{\aleph_1, \aleph_2\}$ such that neither $\aleph_1 \subset \aleph_2$ nor $\aleph_2 \subset \aleph_1$

**Example 7** Consider $\mathfrak{S} = \begin{cases} \aleph_1 = \{(a, 0.10, 0.35, 0.30), (b, 0.20, 0.25, 0.40), (c, 0.50, 0.40, 0.05)\}, \\ \aleph_2 = \{(a, 0.45, 0.30, 0.35), (b, 0.25, 0.30, 0.30), (c, 0.30, 0.35, 0.10)\} \end{cases}$. Then

$$\mathcal{B} = \begin{cases} \aleph_1 = \{(a, 0.10, 0.35, 0.30), (b, 0.20, 0.25, 0.40), (c, 0.50, 0.40, 0.05)\}, \\ \aleph_2 = \{(a, 0.45, 0.30, 0.35), (b, 0.25, 0.30, 0.30), (c, 0.30, 0.35, 0.10)\}, \\ \aleph_1 \cap \aleph_2 = \\ \aleph_3 = \{(a, 0.10, 0.30, 0.35), (b, 0.20, 0.25, 0.40), (c, 0.30, 0.35, 0.10)\} \end{cases}$$ is the base designed from $\mathfrak{S}$,

which further gives

$$\mathfrak{I} = \begin{cases} I = \{(a, 1.00, 0.00, 0.00), (b, 1.00, 0.00, 0.00), (c, 1.00, 0.00, 0.00)\}, \\ O = \{(a, 0.00, 0.00, 1.00), (b, 0.00, 0.00, 1.00), (c, 0.00, 0.00, 1.00)\}, \\ \aleph_1 = \{(a, 0.10, 0.35, 0.30), (b, 0.20, 0.25, 0.40), (c, 0.50, 0.40, 0.05)\}, \\ \aleph_2 = \{(a, 0.45, 0.30, 0.35), (b, 0.25, 0.30, 0.30), (c, 0.30, 0.35, 0.10)\}, \\ \aleph_1 \cap \aleph_2 = \\ \aleph_3 = \{(a, 0.10, 0.30, 0.35), (b, 0.20, 0.25, 0.40), (c, 0.30, 0.35, 0.10)\}, \\ \aleph_1 \cup \aleph_2 = \\ \aleph_4 = \{(a, 0.45, 0.30, 0.30), (b, 0.25, 0.25, 0.30), (c, 0.50, 0.35, 0.05)\}, \\ \aleph_1 \cup \aleph_3 = \\ \aleph_5 = \{(a, 0.10, 0.30, 0.30), (b, 0.20, 0.25, 0.40), (c, 0.50, 0.35, 0.05)\},' \\ \aleph_2 \cup \aleph_3 = \\ \aleph_6 = \{(a, 0.45, 0.30, 0.35), (b, 0.25, 0.25, 0.30), (c, 0.30, 0.35, 0.10)\}, \\ O \cap \aleph_1 = \\ \aleph_7 = \{(a, 0.10, 0.00, 0.30), (b, 0.20, 0.00, 0.40), (c, 0.50, 0.00, 0.05)\}, \\ O \cap \aleph_2 = \\ \aleph_8 = \{(a, 0.45, 0.00, 0.35), (b, 0.25, 0.00, 0.30), (c, 0.30, 0.00, 0.10)\}, \\ O \cap \aleph_3 = \\ \aleph_9 = \{(a, 0.10, 0.00, 0.35), (b, 0.20, 0.00, 0.40), (c, 0.30, 0.00, 0.10)\}, \\ O \cap \aleph_4 = \\ \aleph_{10} = \{(a, 0.45, 0.00, 0.30), (b, 0.25, 0.00, 0.30), (c, 0.50, 0.00, 0.05)\} \end{cases}.$$

5.2 A condition for the smallest picture fuzzy topologies containing $n$ distinct elements

In case of crisp set topology, $\{\emptyset, \aleph_1, \aleph_2, \aleph_3, \ldots, \aleph_{n-2}, X\}$ is the smallest topology containing $\aleph_1, \aleph_2, \aleph_3, \ldots, \aleph_{n-2}$ if $\aleph_1 \subset \aleph_2 \subset \aleph_3 \subset, \ldots, \subset \aleph_{n-2}$. It should be noted that the Pythagorean fuzzy topology follows the same pattern, that is, "$\{O, \aleph_1, \aleph_2, \aleph_3, \ldots, \aleph_{n-2}, I\}$ is the smallest Pythagorean fuzzy topology containing $\aleph_1, \aleph_2, \aleph_3, \ldots, \aleph_{n-2}$ if and only if $\aleph_1 \subset \aleph_2 \subset \aleph_3 \subset, \ldots, \subset \aleph_{n-2}$. It is also true for fuzzy topology and intuitionistic fuzzy topology because Pythagorean fuzzy topology is the generalization of both. But it is easy to verify that it does not holds for picture fuzzy topology of rank greater than one. The following theorem presents a condition for $\{O, \aleph_1, \aleph_2, \aleph_3, \ldots, \aleph_{n-2}, I\}$ to be the smallest picture fuzzy topology containing $\aleph_1, \aleph_2, \aleph_3, \ldots, \aleph_{n-2}$.

**Theorem 13** $\{\emptyset, \aleph_1, \aleph_2, \aleph_3, \ldots, \aleph_{n-2}, X\}$ is the smallest picture fuzyy topology containing $\aleph_1, \aleph_2, \aleph_3, \ldots, \aleph_{n-2}$ if $\aleph_1 \equiv \aleph_2 \equiv \aleph_3 \equiv \cdots \equiv \aleph_{n-2}$.

**Proof.** The proof is an immediate consequence of Theorem 2.

In next example, we show that the converse of Theorem 13 doesn't hold. We construct a picture fuzzy topology of rank 3 with sub-base $\mathfrak{S} = \{\aleph_1, \aleph_2, \aleph_3, \aleph_4\}$ where $\aleph_1 \parallel \aleph_2$, $\aleph_1 \subseteq \aleph_2$, $\aleph_3 \parallel \aleph_4$, $\aleph_3 \subseteq \aleph_4$, $\aleph_2 \nparallel \aleph_3$, $\aleph_1 \equiv \aleph_3$ and $\aleph_2 \equiv \aleph_4$.

**Example 8** Consider $\mathfrak{S} = \begin{cases} \aleph_1 = \{(a, 0.10, 0.15, 0.40), (b, 0.20, 0.10, 0.35), (c, 0.20, 0.15, 0.20)\}, \\ \aleph_2 = \{(a, 0.30, 0.15, 0.35), (b, 0.25, 0.10, 0.30), (c, 0.30, 0.15, 0.10)\}, \\ \aleph_3 = \{(a, 0.10, 0.10, 0.40), (b, 0.20, 0.05, 0.35), (c, 0.20, 0.15, 0.20)\}, \\ \aleph_4 = \{(a, 0.30, 0.10, 0.35), (b, 0.25, 0.05, 0.30), (c, 0.30, 0.15, 0.10)\}, \end{cases}$, then all finite intersections of the elements of $\mathfrak{S}$ evolve

$\mathcal{B} = \begin{cases} \aleph_1 = \{(a, 0.10, 0.15, 0.40), (b, 0.20, 0.10, 0.35), (c, 0.20, 0.15, 0.20)\}, \\ \aleph_2 = \{(a, 0.30, 0.15, 0.35), (b, 0.25, 0.10, 0.30), (c, 0.30, 0.15, 0.10)\}, \\ \aleph_3 = \{(a, 0.10, 0.10, 0.40), (b, 0.20, 0.05, 0.35), (c, 0.20, 0.15, 0.20)\}, \\ \aleph_4 = \{(a, 0.30, 0.10, 0.35), (b, 0.25, 0.05, 0.30), (c, 0.30, 0.15, 0.10)\}, \end{cases}$. The topology constructed from $\mathcal{B}$ is $\mathfrak{I} = \begin{cases} I = \{(a, 1.00, 0.00, 0.00), (b, 1.00, 0.00, 0.00), (c, 1.00, 0.00, 0.00)\}, \\ O = \{(a, 0.00, 0.00, 1.00), (b, 0.00, 0.00, 1.00), (c, 0.00, 0.00, 1.00)\}, \\ \aleph_1 = \{(a, 0.10, 0.15, 0.40), (b, 0.20, 0.10, 0.35), (c, 0.20, 0.15, 0.20)\}, \\ \aleph_2 = \{(a, 0.30, 0.15, 0.35), (b, 0.25, 0.10, 0.30), (c, 0.30, 0.15, 0.10)\}, \\ \aleph_3 = \{(a, 0.10, 0.10, 0.40), (b, 0.20, 0.05, 0.35), (c, 0.20, 0.15, 0.20)\}, \\ \aleph_4 = \{(a, 0.30, 0.10, 0.35), (b, 0.25, 0.05, 0.30), (c, 0.30, 0.15, 0.10)\}, \end{cases}$.

Clearly, $\mathfrak{I} = \{\emptyset, \aleph_1, \aleph_2, \aleph_3, \aleph_4, X\}$ is the smallest picture fuzyy topology containing $\aleph_1, \aleph_2, \aleph_3, \aleph_4$ but $\aleph_1 \equiv \aleph_2 \equiv \aleph_3 \equiv \aleph_4$ is not true.

**Conclusion**

In the present study, we work on picture fuzzy topological spaces. We prove some fundamental results of picture fuzzy subsets and introduce many new notions such as balanced and $\rho$-equivalent relations on picture fuzzy sets. On the basis of these relations, we define rank of picture fuzzy topological spaces to classify picture fuzzy topological spaces. Furthermore, we introduce the concepts of picture fuzzy base and picture fuzzy sub-base. Finally, we develop a technique to form picture fuzzy topological spaces. As far as we know, this is the first endeavor to investigate fuzzy topological spaces. We believe that, the ideas presented in this paper will offer a good platform for further advancement of picture fuzzy topological spaces and their applications in decision making theory.

**Data Availability**

No data were used to support this study.

**Conflicts of Interest**

The authors declare no conflicts of interest regarding the publication of this article.